\documentclass[12pt]{article}
 
\usepackage{babel}
\usepackage{a4wide,doublespace,amstex}
\newtheorem{theorem}{Theorem}[section]
\newtheorem{lemma}[theorem]{Lemma}

\newcommand{\rank}{\mathrm{rank}}
\newcommand{\gl}{\mathrm{gl}}
\newcommand{\vr}{\mathrm{vr}}
\newcommand{\Id}{\mathrm{id}}
\newcommand{\GL}{\mathrm{GL}}

\newtheorem{definition}[theorem]{Definition}
\newtheorem{corollary}[theorem]{Corollary}

\newcommand{\bproof}{\noindent{\bf Proof: }}
\newcommand{\eproof}{\hfill $\Box$\\}
\newcommand{\bremark}{\noindent{\bf Remark: }}
\newcommand{\eremark}{\hfill \\}

\newcommand{\e}{\varepsilon}

\newcommand{\cL}{{\cal L}}

\newcommand{\bN}{{\Bbb N}}
\newcommand{\bR}{{\Bbb R}}

\begin{document}

\title{The Relations between Volume Ratios and\\ New Concepts of $\GL$
Constants}

\author{Y.\ Gordon\footnote{Supported in part by the fund for the
promotion of research in the Technion and by the VPR fund.}
 \and M.\ Junge \and N.J.\ Nielsen\footnote{Supported in part by the Danish
 Natural Science Research Council, grants 9503296 and 9600673.}}

\date{ }
\maketitle

\begin{abstract}
In this paper we investigate a property named $\GL(p,q)$ which is
closely related to the Gordon-Lewis property. Our results on $\GL(p,q)$
are then used to estimate volume ratios relative to $\ell_p$,
$1<p\le\infty$, of unconditional direct sums of Banach spaces.
\end{abstract}

\renewcommand{\theequation}{\arabic{section}.\arabic{equation}}

\section*{Introduction}

In this paper we investigate a property named $\GL(p,q)$, $1\le
p\le\infty$, $1\le q\le\infty$, closely related to the Gordon-Lewis
property $\GL$, and the behavior of $p$-summing norms of operators
defined on direct sums of Banach spaces in the sense of an unconditional
basis. These results are then used to estimate the volume ratios
$\vr(X,\ell_p)$, $1<p\le\infty$, where $X$ is a finite direct sum of
finite dimensional spaces.

A Banach space is said to have $\GL(p,q)$, $1\le p\le\infty$, $1\le
q\le\infty$, if there is a constant $C$ so that $i_q(T)\le C\pi_p(T^*)$
for every finite rank operator $T$ from an arbitrary Banach space to
$X$. Here $\pi_p$ denotes the $p$-summing norm and $i_q$ the
$q$-integral norm. This property was also considered by Reisner
\cite{R}, note however the slight difference in the notation: Our $\GL(p,q)$
corresponds to his $q',p'$-$\GL$-space.

We now wish to discuss the arrangement and contents of this paper in
greater detail.

In Section \ref{sec-1} of the paper we investigate the basic properties
of $\GL(p,q)$ and prove some inequalities for $p$-summing operators,
respectively $q$-integral operators, defined on, respectively with range
in, a direct sum of Banach spaces in the sense of an unconditional
basis. These inequalities are then used to prove that if $(X_n)$ is a
sequence of Banach spaces with uniformly bounded $\GL(p,q)$-constants
and $X$ is the direct sum of the $X_n$'s in the sense of a $p$-convex
and $q$-concave unconditional basis, then $X$ has $\GL(p,q)$ as
well. More generally we obtain that if $Y$ is a Banach space with
$\GL(p,q)$ and $L$ is a $p$-convex and $q$-concave Banach lattice, then
$L(Y)$ has $\GL(p,q)$. $K^p(L)$ and $K_q(L)$ denote the p-convexity and 
q-concavity constants of $L$ respectively.

In Section \ref{sec-2} we combine the results of Section \ref{sec-1}
with those of \cite{GJ} to obtain some estimates of volume ratios. 
One of our results, Theorem \ref{thm-2.5}, has the following geometric consequence:
Let $L$ be a p-convex and q-concave Banach lattice having an $n$-dimensional
Banach space $Y=(\bR^n, \|\cdot\|)$ as an isometric quotient. Let 
$1\le p,q\le\infty$, $\frac{1}{p}+\frac{1}{p'}=1$, then there are $n$-
dimensional linear quotients $V_p$ and $V_q$ of $B_{\ell_p}$ and
$B_{\ell_q}$ respectively, so that $V_q\subseteq B_Y\subseteq V_p$ for which
\[
\Big(\frac{|V_p|}{|V_q|}\Big)^{\frac{1}{n}}\le c\sqrt{p'}\ \gl_{p,q}(L) 
\le c\sqrt{p'}\ K^p(L)K_q(L).
\]

If $X$ is a finite direct sum of $n_k$-dimensional Banach spaces $X_k$,
$1\le k\le m$ in the sense of a finite 1-unconditional basis, then we
prove that 
\[ 
\Big(\prod_{k=1}^m \vr(X_k,\ell_p)^{n_k}\big)^{1/n} \approx
\vr(X,\ell_p)
\]
for $1<p\le\infty$, where $n=\sum_{k=1}^m n_k$.

\setcounter{section}{-1}

\section{Notation and Preliminaries}
\label{sec-0}

In this paper we shall use the notation and terminology commonly used in
Banach space theory as it appears in \cite{LT1}, \cite{LT2}
\cite{P} and \cite{TJ}.

If $X$ and $Y$ are Banach spaces, $B(X,Y)$ ($B(X)=B(X,X)$) denotes the
space of bounded linear operators from $X$ to $Y$ and throughout the
paper we shall identify the tensor product $X\otimes Y$ with the space
of $\omega^*$-continuous finite rank operators from $X^*$ to $Y$ in the
canonical manner. Further if $1\le p<\infty$ we let $\Pi_p(X,Y)$ denote
the space of $p$-summing operators from $X$ to $Y$ equipped with the
$p$-summing norm $\pi_p$, $I_p(X,Y)$ denotes the space of all
$p$-integral operators from $X$ to $Y$ equipped with the $p$-integral
norm $i_p$ and $N_p(X,Y)$ denotes the space of all $p$-nuclear operators
from $X$ to $Y$ equipped with the $p$-nuclear norm $\nu_p$. We recall
that if $1\le p\le\infty$ then an operator $T$ is said to factor through
$L_p$ if it admits a factorization $T=BA$, where $A\in B(X,L_p(\mu))$
and $B\in B(L_p(\mu),Y)$ for some measure $\mu$ and we denote the space
of all operators from $X$ to $Y$, which factor through $L_p$ by
$\Gamma_p(X,Y)$. If $T\in\Gamma_p(X,Y)$ we define
\[
\gamma_p(T) = \inf\{\|A\| \|B\| \mid T=BA, \quad\mbox{$A$ and $B$ as
above}\},
\]
$\gamma_p$ is a norm on $\Gamma_p(X,Y)$ turning it into a Banach
space. All these spaces are operator ideals and we refer to the above
mentioned books and \cite{K}, \cite{PIE} and \cite{KW} for further details. 
To avoid misunderstanding we stress that in this paper a $p$-integral 
operator $T$ from $X$ to $Y$ has a $p$-integral factorization ending in $Y$ 
with $i_p(T)$ defined accordingly; in some books this is referred to as a
strictly $p$-integral operator.

In the formulas below we shall, as is customary, interpret $\pi_\infty$ as 
the operator norm and $i_\infty$ as the $\gamma_\infty$-norm.

If $n\in\bN$ and $T\in B(\ell_2^n,X)$ then following \cite{TJ} we define
the $\ell$-norm of $T$ by
\[
\ell(T) = \big(\int_{\bR^n} \|Tx\|^2 d\gamma(x)\big)^{\frac12}
\]
where $\gamma$ is the canonical Gaussian probability measure on
$\ell_2^n$.

A Banach space $X$ is said to have the Gordon-Lewis property
(abbreviated $\GL$) \cite{GL} if every 1-summing operator from $X$ to an
arbitrary Banach space $Y$ factors through $L_1$. It is easily verified
that $X$ has $\GL$ if and only if there is a constant $K$ so that
$\gamma_1(T)\le K\pi_1(T)$ for every Banach space $Y$ and every $T\in
X^*\otimes Y$. In that case $\GL(X)$ denotes the smallest constant $K$
with this property.

We shall say that $X$ has $\GL_2$ if it has the above property with
$Y=\ell_2$ and we define the constant $\gl(X)$ correspondingly. An easy
trace duality argument yields that $\GL$ and $\GL_2$ are self dual
properties and that $\GL(X)=\GL(X^*)$, $\gl(X)=\gl(X^*)$ when
applicable. It is known \cite{GL} that every Banach space with local
unconditional structure has $\GL$.

If $E$ is a Banach space with a 1-unconditional basis $(e_n)$ and
$(X_n)$ is a sequence of Banach spaces then we put
\[
\big(\sum_{n=1}^\infty X_n\big)_E = \{x\in \prod_{n=1}^\infty X_n\mid
\sum_{n=1}^\infty \|x(n)\|e_n \quad\mbox{converges in $E$}\}
\]
and if $x\in \big(\sum_{n=1}^\infty X_n\big)_E$ we define
\[
\|x\| = \big\|\sum_{n=1}^\infty \|x(n)\| e_n\big\|
\]
thus defining a norm on $(\sum_{n=1}^\infty X_n)_E$ turning it into a
Banach space. If $X_n=X$ for all $n\in\bN$ we put
$E(X)=(\sum_{n=1}^\infty X_n)_E$.

If $(Y_n)$ is another sequence of Banach spaces and $T_n\in B(X_n,Y_n)$
for all $n\in\bN$ with $\sup_n \|T_n\|<\infty$ then we define the
operator $\oplus_{n=1}^\infty T_n: \big(\sum_n X_n\big)_E\to
\big(\sum_{n=1}^\infty Y_n\big)_E$ by \linebreak $\big(\oplus_{n=1}^\infty
T_n\big)(x)=(T_n x(n))$. Clearly $\|\oplus_{n=1}^\infty T_n\|\le \sup_n
\|T_n\|$.

We shall need a ``continuous'' version of the above direct sums so hence
let $X$ be a Banach space and $L$ a Banach lattice. If $\sum_{j=1}^n
x_j\otimes y_j \in X\otimes L$ then it follows from \cite[Section I
d)]{LT2} that $\sup_{\|x^*\|\le 1} \big|\sum_{j=1}^n x^*(x_j)y_j\big|$
exists in $L$ and we put
\[
\big\|\sum_{j=1}^n x_j\otimes y_j\big\|_m = \big\|\sup_{\|x^*\|\le 1}
\big|\sum_{j=1}^n x^*(x_j)y_j\big|\big\|
\]
and define $L(X)$ to be the completion of $X\otimes L$ equipped with the
norm $\|\cdot \|_m$. Spaces of that type was originally defined and
investigated by Schaefer \cite{Sch}; we refer to \cite{HNO} for the
properties of $L(X)$ needed in this paper.

If $n\in\bN$ and $X$ is an $n$-dimensional Banach space then we shall
identify $X$ with $(\bR^n,\|\cdot\|_X)$ by choosing a fixed basis of $X$
and identifying it with the unit vector basis of $\bR^n$, and $B_X$
will denote the unit ball of $X$. Hence if
$B\subseteq X$ is a Borel set we can define the volume $|B|$ of
$B$ as the Lebesgue measure of $B$ considered as a subset of
$\bR^n$. The volume function thus defined is uniquely determined up to a
constant only depending on the chosen basis.

Let $X$ and $Y$ be $n$-dimensional Banach spaces and let
$(x_j)^n_{j=1}$, respectively $(y_j^*)$ be fixed bases of $X$,
respectively $Y^*$. If $T\in B(X,Y)$ then we define the determinant of $T$ by
\[
\det T = \det \{y_j^*(Tx_i)\}.
\]

Up to a constant depending only on the chosen bases $\det T$ is uniquely
determined.

In the sequel, if $X_k$ $1\le k\le m$ are $n_k$-dimensional spaces with
fixed chosen bases and $n=\sum_{n=1}^m n_k$ then we shall always
identify $\prod_{n=1}^m X_k$ with $\bR^n$ via the canonical basis of the
product.

If $X$ is a Banach space and $E$ is an $n$-dimensional Banach space 
then we define the volume ratio $\vr(E,X)$, \cite{GJ}, \cite{GMP} by
\[
\vr(E,X)=\inf\big\{\Big(\frac{|B_E|}{|T(B_X)|}\Big)^{\frac{1}{n}} 
\mid T\in B(X,E), \|T\|\le
1\big\}.
\]
When $X=\ell_\infty$, $\vr(E,\ell_\infty)$ is called the "zonoid" ratio
of $E$, and when $X=\ell_2$, $\vr(E,\ell_2)$ is the well known classical 
volume ratio of $E$, see e.g. \cite{P}, \cite{TJ} and the references therein.

Similarly, 
\[
\vr(E,S(X))=\inf\big\{\Big(\frac{|B_E|}{|T(B_F)|}\Big)^{\frac{1}{n}} 
\mid F\subseteq X, \dim\ F=n,\ T(B_F)\subseteq B_E \big\}.
\]
When $X=\ell_p$ we set $S_p=S(\ell_p)$.

Finally, if $X$ and $Y$ are Banach spaces and $T\in B(X,Y)$ then we
define the $n$-th volume number $v_n(T)$ by
\[
v_n(T) = \sup\big\{\Big(\frac{|T(B_E)|}{|B_F|}\Big)^{\frac{1}{n}}
 \mid E\subseteq X,
T(E)\subseteq F\subseteq Y, \dim E=\dim F=n\big\}.
\]
If $\rank(T) < n$ we put $v_n(T)=0$. Volume numbers or similar notions
were discussed by \cite{D}, \cite{MA}, \cite{PTJ}, \cite{P} 
 and \cite{TJ}. The main results on volume numbers we are
going to use here can be found in \cite{GJ}.

\section{The GL Property and Related Invariances}
\label{sec-1}

We start with the following definition

\begin{definition}
\label{def-1.1}

If $1\le p$, $q\le\infty$ then a Banach space $X$ is said to have
$\GL(p,q)$ if there exists a constant $C$ so that for all Banach spaces
$Z$ and every $T\in Z^*\otimes X$ we have
\begin{equation}
\label{eq-1.1}
i_q(T)\le C\pi_p(T^*).
\end{equation}
\end{definition}

If $X$ has $\GL(p,q)$ then the smallest constant $C$ which can be used
in (\ref{eq-1.1}) is denoted by $\GL_{p,q}(X)$. If $X$ satisfies the
condition of Definition \ref{def-1.1} for $Z=\ell_2$ then we shall say
that $X$ has $\gl(p,q)$ and define the constant $\gl_{p,q}(X)$
correspondingly. It was proved in Corollary (3.12) (I) \cite{GJ} 
that if $X$ is a finite-dimensional Banach space, 
$1\le p,q\le \infty$, and $\frac{1}{p}+\frac{1}{p'}=1$ then
\[ \vr(X,\ell_q)\vr(X^*,\ell_{p'})\le \frac{\pi e}{2}\ \gl_{p,q}(X).
\]
By factoring a given finite rank operator with range in $X$ through its
kernel it is readily seen that it is enough to consider finite
dimensional spaces $Z$ in Definition \ref{def-1.1}. 

By trace duality arguments it is readily verified that if a Banach space
$X$ has $\GL(p,q)$ then $X^*$ has $\GL(q',p')$, and hence $X^{**}$ has
$\GL(p,q)$ as well. The other direction is part of the next lemma.

\begin{lemma}
\label{lemma-1.2}
Let $1\le p,q\le\infty$ and let $X$ be a Banach space. If a subspace $Y$
of $X^{**}$ containing $X$ has $\GL(p,q)$ then $X$ has it as well with
\[
\GL_{p,q}(X)\le \GL_{p,q}(Y).
\]
\end{lemma}

\bproof
Let $Z$ be a finite dimensional Banach space, $T\in Z^*\otimes X$ and $\e>0$ arbitrary. Let
$I$ denote the identity operator of $X$ into $Y$. Choose a finite
dimensional subspace $F$, $IT(Z)\subseteq F\subseteq Y$, 
so that $\e+i_q(IT)\ge i_q(IT: Z\to
F)$. The principle of local reflexivity \cite{JRZ}, \cite{LR}, gives an
isomorphism $V: F\to X$ with $\|V\|\le 1+\e$ and $Vx=x$ for all $x\in F\cap X$.

Since $VIT=T$ we obtain
\begin{eqnarray*}
i_q(T) &\le & \|V\| i_q(IT: Z\to F)\\
&\le & (1+\e)i_q(IT)+(1+\e)\e\\
&\le & (1+\e)GL_{p,q}(Y)\pi_p(T^*I^*)+(1+\e)\e\\
&\le & (1+\e)\pi_p(T^*)+(1+\e)\e.
\end{eqnarray*}
Since $\e$ was arbitrary this shows that $X$ has $\GL(p,q)$ with
$\GL_{p,q}(X)\le \GL_{p,q}(Y)$.
\eproof

It follows
immediately that $X$ has $\GL(\infty,q)$ for some $q$, $1\le q < \infty$
 (or dually has $\GL(p,1)$ for some $p$, $1\le p<\infty$) if and only 
 if it is finite dimensional. Obviously, since the $\GL$-property is self-dual, $X$ has $\GL$
if and only if it has $\GL(1,\infty)$.

The next theorem which is the result of the work of several authors,
\cite{JRZ}, \cite{KW}, \cite{LR} and \cite{MP}, describes the situation 
for the remaining values of $p$ and $q$.

\begin{theorem}
\label{thm-1.2}
If $X$ is a Banach space, $1\le p,q\le\infty$, then the following
statements hold:
\begin{itemize}
\item[(i)] If $X$ has $\GL(p,q)$ then $X$ has $\GL$, $X$ is of cotype
$\max(q,2)$ and $X^*$ is of cotype $\max(p',2)$. If $q<\infty$ and
$1<p<\infty$, then $X$ is of type $\min(2,p)$.
\item[(ii)] If $X$ has $\GL$, $2\le q<\infty$ and
$B(L_\infty,X)=\Pi_q(L_\infty,X)$ then $X$ has $\GL(1,q)$.
\item[(iii)] If $X$ has $\GL$ and is of type $p$-stable for some $p$,
$1<p\le 2$, then there is a $q$, $1\le q<\infty$ so that $X$ has
$\GL(p,q)$.
\item[(iv)] If $1<p<\infty$ then $X$ has $\GL(p,p)$ if and only if $X$
is either a $\cL_p$-space or isomorphic to a Hilbert space.
\item[(v)] If $1<q<p<\infty$ then $X$ has $\GL(p,q)$ if and only if
$X$ is isomorphic to a Hilbert space.
\item[(vi)] $X$ has $\GL(\infty,\infty)$ (respectively $X$ has $\GL(1,1)$)
 if and only if it is a $\cL_\infty$-space (respectively a $\cL_1$-space).
\item[(vii)] If $X$ is a $p$-convex and $q$-concave Banach lattice then
$X$ has $\GL(p,q)$ with $\GL_{p,q}(X)\le K^p(X)K_q(X)$.

\end{itemize}
\end{theorem}
\bproof
\begin{itemize} 
\item[(i)]  Assume that $X$ has $\GL(p,q)$. If $T\in Y^*\otimes X$
 then
\[ 
 \gamma_1(T^*)=\gamma_\infty (T)\le i_q(T)\le \GL_{p,q}(X)\pi_p(T^*)
 \le \GL_{p,q}(X)\pi_1(T^*)
 \]
 and hence $X^*$ and therefore also $X$ has $\GL$.
 
   If $q$ is finite then $X$ has
property $(S_q)$ of \cite{CN} and is therefore of cotype $\max(q,2)$ by
Theorem 1.3 there. Since $X^*$ has $\GL(q',p')$ it is of cotype
$\max(p',2)$.

If $q<\infty$ and $1<p<\infty$, then both $X$ and $X^*$ are of finite
cotype and since in addition $X$ has $\GL$ it follows e.g.\ from
\cite[Theorem 1.9]{CN} that $X$ is $K$-convex. Therefore $X$ has type
$\min(p,2)$.
\item[(ii)] Assume that $X$ has $\GL$ and
$\Pi_q(L_\infty,X)=B(L_\infty,X)$ for some $q$, $2\le q<\infty$ with
$K$-equivalence between the norms and let $C$ be the $\GL$-constant of
$X$.

Let $Z$ be an arbitrary finite dimensional Banach space and $T\in Z^*\otimes X$. Since $X$
and hence also $X^*$ has $\GL$ with constant $C$, there exists a measure
$\mu$ and operators $A\in B(X^*,L_1(\mu))$, $B\in B(L_1(\mu),Z^*)$ with
$T^*=BA$ and $\|A\| \|B\|\le C\pi_1(T^*)$. 

Let $\e >0$ be arbitrary. Using the local properties of
$L_\infty(\mu)$ we can find a finite dimensional subspace $E\subseteq
L_\infty(\mu)$ with $d(E,\ell_\infty^{\dim E})\le 1+\e$ and
$B^*(Z)\subseteq E$, and hence $i_q(A^*_{|E})\le
(1+\e)\pi_q(A^*_{|E})$. By the principle of local reflexivity
there is an isomorphism $V: A^*(E)\to X$ so that $\|V\|\le 1+\e$ and
$Vx=x$ for all $x\in A^*(E)\cap X$. Since clearly $T=VA^*_{|E}B^*$ we
obtain
\[
i_q(T)\le (1+\e)^2\|B\| \pi_q(A^*_{|E})\le K(1+\e)^2 \|A\| \|B\|\le
KC(1+\e)^2 \pi_1(T^*)
\]
and hence $X$ has $\GL(1,q)$.
\item[(iii)] Let $X$ have $\GL$ and be of type $p$-stable for some $p$,
$1<p\le 2$. By \cite{MP} $\pi_{p'}(L_\infty,X^*)=B(L_\infty,X^*)$ and
$X$ is of finite cotype so that there is a finite $q$ with
$B(L_\infty,X)=\pi_q(L_\infty,X)$. The first statement implies that
$\Pi_1(X^*,Z) = \Pi_p(X^*,Z)$ for any Banach space $Z$ and therefore $X$
has $\GL(p,q)$ by (ii).
\item[(iv)] Let $1<p<\infty$. By \cite{KW} $X$ has $\GL(p,p)$ if and
only if $X$ is isomorphic to a complemented subspace of an $L_p$-space,
or equivalently \cite{LR} if and only if either $X$ is a $\cL_p$-space
or isomorphic to a Hilbert space.
\item[(v)] Let $1<q<p<\infty$. Since $q>1$ it follows from
\cite[Proposition 0.3]{CN} and its proof that there is a universal 
constant $c$ so that if $Z$ is a Banach space and $T\in Z^*\otimes \ell_2$
then $\frac{1}{c\sqrt q} i_q(T)\le \ell(T^*)\le c\sqrt p\ \pi_p(T^*)$. 
This gives that $\ell_2$ has $\GL(p,q)$.

If $X$ has $\GL(p,q)$ then it has both $\GL(p,p)$ and $\GL(q,q)$ and
hence it follows from (iv) that $X$ is isomorphic to a Hilbert space.
 \item[(vi)] 
 Assume that $X$ is a $\cL_{\infty,\lambda}$-space and let
$T\in Z^*\otimes X$, where $Z$ is an arbitrary Banach space.
By definition, there is a finite dimensional subspace $E\subseteq X$ with
$d(E,\ell_\infty^{\dim E})\le \lambda$ and $T(Z)\subseteq E$. Hence $\gamma_\infty(T)\le 
\lambda ||T||$ and $X$ has $\GL(\infty,\infty)$.

Assume next that $X$ has $\GL(\infty,\infty)$. By using trace duality twice, 
we obtain that the identity operator of $X^{**}$ factors through an $L_\infty$-
space and is therefore isomorphic to a complemented subspace of an 
$L_\infty$-space. It follows from \cite{LR} that $X^{**}$ and hence $X$
is $\cL_\infty$-space.

By Lemma \ref{lemma-1.2} $X$ has $\GL(1,1)$ if and only if $X^*$ has $\GL(\infty,\infty)$
and the result follows by noting that $X$ is a $\cL_1$-space if and only
if $X^*$ is a $\cL_\infty$-space.

\item[(vii)] Let $Y$ be an arbitrary Banach space and $T\in Y^*\otimes X$.
 From \cite[Theorem 1.3]{HNO} or \cite{J} it follows that 
 there exists a factorization
 $T=BA$, where $A\in B(Y,L_\infty)$, $B\in B(L_\infty,X)$, $B\ge 0$ and
 $\|A\| \|B\|\le K^p(X)\pi_p(T^*)$. By a theorem of Maurey, see e.g.
 \cite[Theorem 1.d.10]{LT2}, $B$ is $q$-summing, with $i_q(B)=
 \pi_q(B)\le K_q(X) \|B\|$. Hence, $i_q(T)\le \|A\|i_q(B)\le
 K^p(X)K_q(X)\pi_p(T^*)$. The case $p=1,\ q=\infty$ was first proved in
 \cite{GL}, the general formula is contained in \cite{J} up to trace 
 duality.
\end{itemize}
\eproof

In the rest of this section we let $E$ denote a Banach space with a normalized
1-unconditional basis $(e_j)$ and biorthogonal system $(e_j^*)$. Our first
result is a generalization of \cite[Theorem 1.18]{CN}.

\begin{theorem}
\label{thm-1.3}
Let $(X_n)$ be a sequence of Banach spaces, $Z$ a Banach space, $1\le
p\le\infty$ and $T: \left(\sum_{n=1}^\infty X_n\right)_E\to Z$ a
$p$-summing operator. Put $T_n = T_{|X_n}$ for all $n\in\bN$.

If $E^*$ is $p$-convex then for all $m\in \bN$
\begin{equation}
\label{eq-1.2}
\big\|\sum_{j=1}^m \pi_p(T_j)e^*_j\big\|\le K^{p}(E^*)\pi_p(T).
\end{equation}

If, in addition, $(e^*_j)$ is boundedly complete or $\sum^\infty_{n=1}
T_n$ converges to $T$ in the $p$-summing norm, then $\sum^\infty_{j=1}
\pi_p(T_j)e^*_j$ converges in $E^*$ and
\begin{equation}
\label{eq-1.3}
\big\|\sum_{j=1}^\infty \pi_p(T_j)e^*_j\big\|\le K_{p'}(E) \pi_p(T).
\end{equation}
\end{theorem} 

\bproof
It is obvious that (\ref{eq-1.2}) implies (\ref{eq-1.3}) under the
additional assumptions, so let us concentrate on (\ref{eq-1.2}).

Let $\e>0$ be given arbitrarily. For every $n\in\bN$ we can choose a
finite set $\sigma_n\subseteq\bN$ and $\{x_i(n)\mid
i\in\sigma_n\}\subseteq X_n$ so that
\begin{gather}
\label{eq-1.4}
\pi_p(T_n)^p\le\sum_{i\in\sigma_n} \|Tx_i(n)\|^p+\e\cdot 2^{-n},
\label{eq-1.5} \\
\sup\big\{\sum_{i\in\sigma_n}|x^*(x_i(n))|^p \mid x^*\in X^*_n,
\|x^*\|\le 1\big\} \le 1.
\end{gather}
For every sequence $(\alpha_n)\subseteq\bR_+\cup\{0\}$ and every $m\in\bN$ 
we obtain:
\begin{multline}
\label{eq-1.6}
\sum^m_{n=1} \alpha_n\pi_p(T_n)^p \le  \sum_{n=1}^m \sum_{i\in\sigma_n}
\|T(\alpha_n^{1/p} x_i(n)\|^p+\e \\
\le  \pi_p(T)^p \sup\big\{\sum^m_{n=1} \alpha_n \sum_{i\in\sigma_n}
|\langle x^*(n),x_i(n)\rangle|^p \mid x^*  \in 
\big(\sum_{n=1}^{\infty} X_n^*\big)_{E^*},\|x^*\|\le 1\big\} +\e \\
 \le  \pi_p(T)^p \sup\big\{ \sum^m_{n=1} \|x^*(n)\|^p \alpha_n
\sum_{i\in\sigma_n} \mid\big\langle
\frac{x^*(n)}{\|x^*(n)\|},x_i(n)\big\rangle|^p\mid x^*  \in
\big(\sum^\infty_{n=1}X_n^*\big)_{E^*}, \|x^*\|\le 1\big\}+\e
\\
\le  \pi_p(T)^p \sup \big\{\sum_{n=1}^m \|x^*(n)\|^p \alpha_n \mid
x^* \in \big(\sum_{n=1}^\infty X_n^*\big)_{E^*},\|x^*\|\le
1\big\}+\e.
\end{multline}
Let $E^*_{(p)}$ denote the $p$-concavification of $E^*$ (see \cite{LT2}
for details). If in (\ref{eq-1.6}) we take the supremum over all
sequences $(\alpha_n)\in (E^*_{(p)})^*$ with $\|(\alpha_n)\|\le 1$ and
let $\e\to 0$ we get
\begin{equation}
\label{eq-1.7}
\|\sum^m_{n=1} \pi_p(T_n)e_n^*\|^p \le (K^p(E^*))^p\pi_p(T)^p
\end{equation}
which is (\ref{eq-1.2}).
\eproof

By the trace duality between the $p'$-integral and the $p$-summing norms
we obtain with the same notation as in Theorem \ref{thm-1.3}:

\begin{corollary}
\label{cor-1.3}
Let $E$ be $q$-concave for some $q$, $1\le q\le\infty$, $S\in
B(Z,(\sum_{j=1}^\infty X_j)_E)$ and denote by $P_n: (\sum_{j=1}^\infty
X_j)_E\to X_n$ the canonical projection for all $n\in\bN$.

If $\sum_{n=1}^{\infty}i_q(P_nS)e_n$ converges in $E$ then 

\begin{equation}
\label{eq-1.8}
i_q(S) \le K_q(E) \|\sum_{n=1}^\infty i_q(P_nS)e_n\|_E.
\end{equation}
\end{corollary}

The analogous result holds for $q$-nuclear operators.

\bproof
Let $\e>0$ be arbitrary. For every $n\in\bN$ we can find a measure
$\mu_n$, $A_n\in B(Z,L_q(\mu_n))$ and $B_n\in B(L_q(\mu_n),X_n)$ so that
$i_q(A_n)\le i_q(P_nS)+\e\cdot 2^{-n}$, $\|B_n\|\le 1$ and
$P_nS=B_nA_n$.

Let $m\in\bN$ and consider the $q$-integral operator $\sum^m_{n=1} A_n:
Z\to (\sum^m_{n=1}L_q(\mu_n))_E$. By trace duality \cite{PP} we can find
a $q'$-summing operator $W: (\sum^m_{n=1} L_q(\mu_n))_E\to Z$ with
$\pi_{q'}(W)=1$ so that with $W_n = W_{|L_q(\mu_n)}$ we have 
\begin{eqnarray}
\label{eq-1.9}
i_q(\sum^m_{n=1} A_n) &\le & \sum_{n=1}^m tr(W A_n)+\e\le \sum^m_{n=1}
\pi_{q'}(W_n)i_q(A_n) +\e \nonumber \\
&\le & \|\sum^m_{n=1} \pi_{q'}(W_n)e^*_n\| (\|\sum^m_{n=1}
i_q(P_nS)e_n\|+\e)+\e \nonumber \\
&\le & K_q(E)\|\sum^m_{n=1} i_q(P_nS)e_n\|+\e(K_q(E)+1)
\end{eqnarray}
where we have used Theorem \ref{thm-1.3} to get the last inequality.

Formula (\ref{eq-1.9}) shows that $\sum^\infty_{n=1} A_n$ converges in the
$q$-integral norm to an operator $A: Z\to (\sum^\infty_{n=1}
L_q(\mu))_E$ with
\begin{equation}
\label{eq-1.10}
i_q(A)\le K_q(E) \|\sum^\infty_{n=1} i_q(P_nS)e_n\|+\e(1+K_q(E)).
\end{equation}
The operator $B = \oplus^\infty_{n=1} B_n:
(\sum^\infty_{n=1}L_q(\mu_n))_E\to (\sum^\infty_{n=1} X_n)_E$ is clearly
bounded with $\|B\|\le 1$ and $S=BA$. Hence $S$ is $q$-integral with
\begin{equation}
\label{eq-11}
i_q(S)\le K_q(E) \|\sum^\infty_{n=1} i_q(P_nS)e_n\|.
\end{equation}
The statement on $q$-nuclear operators can be proved in a similar
manner or by duality.
\eproof

The next theorem was originally proved by Reisner \cite{R} using other
methods; we can also obtain it directly from Theorem
\ref{thm-1.2} and Corollary \ref{cor-1.3}.

\begin{theorem}
\label{thm-1.5}
Let $1\le p\le q\le\infty$ and let $(X_n)$ be a sequence of Banach
spaces all having $\GL(p,q)$ with $M=\sup\{\GL_{p,q}(X_n)\mid
n\in\bN\}<\infty$. If $E$ is $p$-convex and $q$-concave then
$(\sum^\infty_{n=1}X_n)_E$ has $\GL(p,q)$ with
\begin{equation}
\label{eq-1.12}
\GL_{p,q}((\sum^\infty_{n=1}X_n)_E)\le MK^p(E)K_q(E).
\end{equation}
\end{theorem}
\bproof
Let $Z$ be an arbitrary Banach space and $T\in Z^*\otimes
(\sum^\infty_{n=1}X_n)_E$. Put
$S=T^*_{|(\sum^\infty_{n=1}X_n^*)_{E^*}}$. From Theorem \ref{thm-1.3}
and Corollary \ref{cor-1.3} we obtain
\begin{eqnarray}
\label{eq-1.13}
i_q(T) &\le & K_q(E)\|\sum^\infty_{n=1} i_q(P_nT)e_n\| \nonumber \\
&\le & M K_q(E)\|\sum^\infty_{n=1}\pi_p(T^*P^*_n)e_n\|\le MK^p
(E)K_q(E)\pi_p(S) \nonumber \\
&\le & MK^p(E)K_q(E)\pi_p(T^*),
\end{eqnarray}
from which the statement directly follows.
\eproof

The case $p=1$ and $q=\infty$ gives:

\begin{corollary}
\label{cor-1.5}
If $(X_n)$ is a sequence of Banach spaces all having $\GL$ with
$M=\sup\{\GL(X_n)\mid n\in\bN\}<\infty$ then $(\sum^\infty_{n=1}X_n)_E$
has $\GL$ with
\begin{equation}
\label{eq-1.14}
\GL((\sum^\infty_{n=1}X_n)_E)\le M.
\end{equation}
\end{corollary}

We now wish to generalize Theorem \ref{thm-1.5} and its corollary to the
space $L(X)$, where $L$ is a Banach lattice and $X$ a Banach space. For
this we need the following lemma on $p$-summing norms which might
also be useful in other situations. Before we state it we need a little
notation: If $X$ and $Y$ are Banach spaces, $T\in X^*\otimes Y$ and $F$
is a finite dimensional subspace of $Y$ with $T(X)\subseteq F$ then
$T_F$ denote the operator $T$ considered as an operator from $X$ to $F$.

\begin{lemma}
\label{lemma-1.7}
Let $X$ and $Y$ be Banach spaces and $T\in X^*\otimes Y$. If $1\le
p\le\infty$ and ${\cal F}$ is an upwards directed set of finite
dimensional subspaces all containing $T(X)$ with $\overline{\cup\{F\mid
F\in{\cal F}\}}=Y$, then
\begin{equation}
\label{eq-1.15}
\lim_{F\in{\cal F}} \pi_p(T^*_F) = \pi_p(T^*).
\end{equation}
\end{lemma}

\bproof
We can without loss of generality assume that $X$ is finite dimensional
and let us also assume that $1\le p<\infty$; the case $p=\infty$ is
easier and left to the reader.

The net $(\pi_p(T^*_F))$ is non-negative and decreasing and hence
convergent to $\alpha$, say; clearly $\pi_p(T^*)\le\alpha$.

Let now $\e>0$ be arbitrary. Since $\rank(T^*_F)=\rank(T^*)$ for all $F$
it follows from \cite[Theorem 5 and its proof]{DJ} that we can find an 
$m$ independent
of $F$, so that $\pi_p(T^*_F)$ can be computed up to $\e$ using $m$
vectors from $F^*$. Hence for every $F\in{\cal F}$ we can find
$\{x^*_{j,F}\mid 1\le j\le m\}\subseteq X^*$ so that
\begin{gather}
\label{eq-1.16}
\|x^*_{j,F}\|  <  1 \quad\mbox{for all} \quad 1\le j\le m\\
\label{eq-1.17}
\sum^m_{j=1} |x^*_{j,F}(x)|^p < 1\quad\mbox{for all}\quad x\in F,\quad
\|x\|\le 1\\
\label{eq-1.18}
\left(\sum^m_{j=1} \|T^*x^*_{j,F}\|^p\right)^{1/p} \ge  \pi_p(T^*_F)-\e.
\end{gather}
(\ref{eq-1.16}) gives that there is a subnet $(x^*_{j,F'})$ and an
$x^*_j\in X^*$ so that $(x^*_{j,F'})$ converges $w^*$ to $x^*_j$ for all
$j$, $1\le j\le m$. Since $\overline{\cup_{F\in{\cal F}} F}=Y$,
(\ref{eq-1.17}) gives that $\sum^m_{j=1} |x_j^*(x)|^p\le 1$ for all $x\in
Y$, $\|x\|\le 1$.

From the $w^*$-continuity of $T^*$ it follows that $(T^*x^*_{j,F'})$
converges $w^*$ to $T^*x^*_j$ and therefore also in norm, since $X$ is
finite dimensional. Hence going to the limit in (\ref{eq-1.18}) we get
\[
\pi_p(T^*)\ge\left(\sum^m_{j=1} \|T^*x^*_j\|^p\right)^{1/p}\ge \alpha-\e
\]
which implies that $\pi_p(T^*)\ge\alpha$, since $\e$ was arbitrary.
\eproof

We are now able to prove:

\begin{theorem}
\label{thm-1.8}
Let $1\le p\le q\le \infty$ and let $X$ be a Banach space with
$\GL(p,q)$. If $L$ is a $p$-convex and $q$-concave Banach lattice then
$L(X)$ has $\GL(p,q)$ with
\begin{equation}
\label{eq-1.19}
\GL_{p,q}(L(X))\le\GL_{p,q}(X)K^p(L)K_q(L)
\end{equation}
\end{theorem}

\bproof
If the statement of the theorem has been proved for order complete
Banach lattices, then since $L^{**}$ is order complete and
$L(X)\subseteq L^{**}(X)\subseteq L(X)^{**}$ it follows from Lemma
\ref{lemma-1.2} that $L(X)$ has $\GL(p,q)$. It is therefore no
restriction to assume that $L$ is order complete.

Let $Z$ be a Banach space, $T\in Z^*\otimes L(X)$ with $\|T\|\le 1$ and
$\e>0$ arbitrary. From Lemma \ref{lemma-1.7} it follows that there is
an $n\in{\mathbb N}$ and an $n$ dimensional subspace $F\subseteq L(X)$
so that $T(Z)\subseteq F$ and
\begin{equation}
\label{eq-1.20}
\pi_p(T^*_F)\le \pi_p(T^*)+\e.
\end{equation}
Let $(u_j)^n_{j=1}$ be an Auerbach basis, \cite{LT1}, of $F$ with
biorthogonal system $(u^*_j)\subseteq F^*$. By \cite[Lemma 2.15]{HNO},
and the order completeness of $L$ 
there is an $m\in\bN$, a set $\{e_i\mid 1\le i\le m\}\subseteq L$
consisting of mutually disjoint positive vectors of norm 1 
and $\{v_j\mid 1\le
j\le n\}\subseteq [e_i](X)=Y$ (naturally considered as a subspace of
$L(X)$) so that
\begin{equation}
\label{eq-1.21}
\|u_j-v_j\|\le\frac{\e}{n}\quad\mbox{for all}\quad 1\le j\le n.
\end{equation}
If $S=\sum^n_{j=1} u^*_j\otimes v_j: F\to Y$, then for all $u\in F$ we
have
\begin{equation}
\label{eq-1.22}
\|u-Su\|\le\sum^n_{j=1} \|u_j^*\| \|u_j-v_j\|\le \e.
\end{equation}
Considering $T-ST$ as an operator from $Z$ to $L(X)$ it has the
representation
\begin{equation}
\label{eq-1.23}
T-ST = \sum^n_{j=1} T^*_F u^*_j\otimes (u_j-v_j)
\end{equation}
and therefore
\begin{equation}
\label{eq-1.24}
\nu_1(T-ST) \le \sum^n_{j=1} \|T^*_F u^*_j\| \|u_j-v_j\|\le\e.
\end{equation}
Applying Theorem \ref{thm-1.5} and the previous inequalities we obtain
\begin{eqnarray}
\label{eq-1.25}
i_q(T) &\le & i_q(ST)+i_q(T-ST) \le i_q((ST)_Y)+\nu_1(T-ST)\nonumber \\
&\le & \GL_{p,q}(X)K^p(L)K_q(L)\pi_p((ST)^*_Y)+\e \nonumber \\
&\le & \GL_{p,q}(X)K^p(L)K_q(L)\|S\|\pi_p(T_F^*)+\e \nonumber \\
&\le & \GL_{p,q}(X)K^p(L)K_q(L)(1+\e)(\pi_p(T^*)+\e)+\e.
\end{eqnarray}
From (\ref{eq-1.25}) we conclude that $L(X)$ has $\GL(p,q)$ and since
$\e$ was arbitrary (\ref{eq-1.19}) follows.
\eproof

As a corollary we obtain

\begin{corollary}
\label{cor-1.9}
Let $X$ be a Banach space with $\GL$ and $L$ a Banach lattice. 
Then $L(X)$ has $\GL$ with
\[
GL(L(X)) = \GL(X).
\]
\end{corollary}

\section{Volume Ratios of Direct Sums of Banach Spaces}
\label{sec-2}
\setcounter{equation}{0}

In this section we shall use the results of Section \ref{sec-1} to
compute volume ratios of certain direct sums of finite dimensional
Banach spaces. Throughout the section we let $m\in\bN$, $n_k\in\bN$ for
all $1\le k\le m$, let $E$ denote an $m$ dimensional Banach space with a normalized
1-unconditional basis $(e_k)^m_{k=1}$ and biorthogonal system
$(e_k^*)^m_{k=1}$ and let $X_k$ and $Y_k$ be Banach spaces with $\dim
X_k=\dim Y_k=n_k$ for all $1\le k\le m$; put $n=\sum^m_{k=1}n_k$,
$X=\big(\sum^m_{k=1}X_k\big)_E$ and $Y=\big(\sum^m_{k=1} Y_k\big)_E$. We
wish to compute $\vr(X,\ell_p)$ for $1< p\le\infty$, but before we
can do that we need a few lemmas.

\begin{lemma}
\label{lemma-2.1}
\[
\frac{|B_X|}{|B_Y|} = \prod^m_{k=1}\frac{|B_{X_k}|}{|B_{Y_k}|}
\]
\end{lemma}

\bproof
We will iteratively interchange the $X_k$'s by the $Y_k$'s. 
So, we define $Z_{m-1}=(\sum_{k=1}^{m-1} X_k \oplus Y_m)_E$.
By choosing a basis in each of the involved spaces we may identify
$X_k,Y_k$ respectively $X,Y,Z_{m-1}$ with $\bR^{n_k}$ respectively $\bR^n$ in a
canonical manner. Then  we define $r: \bR^{m-1}\to\bR$ by
\begin{equation}
\label{eq-2.2}
r(t_1,t_2,\dots,t_{m-1}) = \inf\{t\in\bR\mid \|\sum^{m-1}_{j=1}
t_je_j+te_m\|_X=1\}\quad\mbox{for all}\
(t_1,t_2,\dots,t_{m-1})\in\bR^{m-1}.
\end{equation}
For every $x_k\in X_k$, $1\le k\le m-1$ we put
\begin{equation}
\label{eq-2.3}
A(x_1,x_2,\dots,x_{m-1})=r(\|x_1\|_{X_1},\|x_2\|_{X_2},\dots,
\|x_{m-1}\|_{X_{m-1}})B_{X_m}
\end{equation}
and consider $\bR^n=\prod^m_{k=1}\bR^{n_k}$. With this notation we
obtain
\begin{eqnarray}
\label{eq-2.4}
|B_X| &=& \int_{\bR^{n_m}}\dots \int_{\bR^{n_1}}{\bf 1}_{B_X}(x_1,
\dots,x_m)dx_1\dots dx_m
\nonumber \\
&=& \int_{\bR^{n_{m-1}}}\dots \int_{\bR^{n_1}}\Big[\int_{\bR^{n_m}} 
{\bf 1}_{A(x_1,\dots,x_{m-1})}(x_m)
dx_m\Big] dx_1 \dots dx_{m-1} 
\nonumber \\
&=&|B_{X_m}|\int_{\bR^{n_{m-1}}} \dots \int_{\bR^{n_1}}r(\|x_1\|_{X_1},
\dots,\|x_{m-1}\|_{X_{m-1}})
dx_1 \dots dx_{m-1}
\end{eqnarray}
Using the same calculation for
$Z_{m-1}$ 
 in (\ref{eq-2.4}) we get
\[
|B_{Z_{m-1}}| = |B_{Y_m}| \int_{\bR^{n_{m-1}}}\dots \int_{\bR^{n_1}}
r(\|x_1\|_{X_1},\dots,\|x_{m-1}\|_{X_{m-1}})dx_1\dots dx_{m-1}
\]
and hence:
\begin{equation}
\label{eq-2.5}
|B_X| = |B_{Z_{m-1}}|\frac{|B_{X_m}|}{|B_{Y_m}|}.
\end{equation}
The result now follows by iterating (\ref{eq-2.5}).
\eproof

\begin{lemma}
\label{lemma-2.2}
Let $1\le p<\infty$ and for every $1\le k\le m$ let $T_k\in
X_k^*\otimes\ell_2$. If $T=\bigoplus^m_{k=1}T_k\in X^*\otimes
\ell_2^m(\ell_2)$ then
\begin{equation}
\label{eq-2.6}
\pi_p(T) \le c\sqrt{p}\ \big\|\sum^m_{k=1}\pi_p(T_k)e_k^*\big\|
\le c\sqrt{p}\ K_{p'}(E)\pi_p(T)
\end{equation}
where $c$ is a universal constant.
\end{lemma}

\bproof
The right hand side inequality (2.5) is formula (1.2). 
If $2<p$ then it follows from Maurey's extension theorem \cite{M}
that every $p$-summing
operator from an arbitrary Banach space to $\ell_2$ is already
2-summing and $\pi_2(T)\le c\sqrt p\ \pi_p(T)$ (see e.g. \cite{P2}). 
Therefore it suffices to prove (2.5) for $1\le p\le 2$.

By the factorization theorem of Pietsch (see e.g.\ \cite{LT1}) we can
for every $1\le k\le m$ find a Radon probability measure $\mu_k$ on
$B_{X_k^*}$ so that for every $x_k\in X_k$,
\[
\|T_kx_k\|\le \pi_p(T_k)\big(\int_{B_{X_k^*}}|x_k^*(x_k)|^p
d\mu_k(x_k^*)\big)^{\frac{1}{p}}.
\]
Put $B=\prod^m_{k=1} B_{X_k^*}$, $\mu=\prod^m_{k=1}\mu_k$ and
$\tau=\big\|\sum^m_{k=1}\pi_p(T_k)e_k^*\big\|$ and let $(r_k)$
denote the sequence of Rademacher functions on $[0,1]$. Since
$(e_k^*)$ is 1-unconditional the function $f$ defined by
 $f(t,x^*)=\tau^{-1}(r_k(t)\pi_p(T_k)x_k^*)^m_{k=1}$ 
 for all $t\in [0,1]$ and all
$x^*=(x_k^*)\in X^*$ maps $[0,1]\times B$ into $B_{X^*}$ and hence
$d\nu=[dt\times d\mu]\circ f^{-1}$ is a probability measure 
on $B_{X^*}$ (actually concentrated on the sphere).

Since the cotype 2 constant of $L_p(\mu)$ is majorized by $\sqrt 2$ for $1\le p\le 2$ (the constant in Khintchine's inequality for $p=1$, see \cite{H} and \cite{Sz}),
the following inequalities hold for every $x=(x_k)\in X$ (putting
$x^*=(x_k^*)\in X^*$):
\begin{eqnarray}
\|Tx\| &=& \big(\sum^m_{k=1} \|T_kx_k\|^2\big)^{\frac12} \le
\big[\sum^m_{k=1} \pi_p(T_k)^2\big(\int_B |x_k^*(x_k)|^p
d\mu(x^*)\big)^{\frac{2}{p}}\big]^{\frac12} \nonumber \\
&\le & \sqrt 2\ \big(\int_0^1\int_B \big|\sum^m_{k=1}
r_k(t)\pi_k(T_k)x_k^*(x_k)\big|^p d\mu(x^*)dt\big)^{\frac{1}{p}} 
\nonumber \\
&=& \sqrt 2\ \tau\big(\int_{B_{X^*}} |x^*(x)|^p d\nu(x^*)\big)^{\frac{1}{p}}.
\end{eqnarray}
This shows that $T$ is $p$-summing with $\pi_p(T)\le c\sqrt p\ \tau$ for all
$1\le p<\infty$.
\eproof

\begin{lemma}
\label{lemma-2.3}
For every $1\le k\le m$ there exist $\alpha_k\ge 0$, $\beta_k\ge 0$ so
that 
\[
\big\|\sum^m_{k=1}\alpha_ke^*_k\big\|=1=\big\|\sum^m_{k=1}
\beta_ke_k\big\|\]
and 
\[\prod^m_{k=1}
(\alpha_k^{n_k}\beta_k^{n_k} n_k^{-n_k})\ge n^{-n} . \]
\end{lemma}

\bproof
Consider the space $Z = \big(\sum^m_{k=1} \ell_1^{n_k}\big)_{E^*}$ with
its canonical 1-unconditional basis. Applying a result of Lozanovskii
\cite{L} on 1-unconditional bases (see Corollary 3.4 in \cite{P}) 
to $Z$ we see
that there exist numbers $\tau_{jk}\ge 0$ and $\sigma_{jk}\ge 0$ for
$1\le k\le m$ and $1\le j\le n_k$ so that if $D_\tau: \ell^n_\infty \to
Z$, respectively $D_\sigma: Z\to \ell_1^n$ are the diagonal operators
determined by $(\tau_{jk})$, respectively $(\sigma_{jk})$, then:
\begin{equation}
\label{eq-2.9}
\|D_\tau\|=1=\|D_\sigma\|;\qquad D_\tau D_\sigma=\frac{1}{n} \Id_Z.
\end{equation}
If we for every $1\le k\le m$ define
\begin{equation}
\label{eq-2.10}
\alpha_k = \sum^{n_k}_{j=1} \tau_{jk};\qquad \beta_k=\max_{1\le j\le
n_k} \sigma_{jk},
\end{equation}
then
\begin{equation}
\label{eq-2.11}
\big\|\sum^m_{k=1} \alpha_k e_k^*\big\| = \|D_\tau\|=1=\|D_\sigma\|=
\big\|\sum^m_{k=1} \beta_k e_k\big\|.
\end{equation}
From (\ref{eq-2.9}) we obtain for all $1\le k\le m$:
\begin{equation}
\label{eq-2.12}
\frac{1}{n} = \big(\prod^k_{j=1}
\sigma_{jk}\tau_{jk}\big)^{\frac{1}{n_k}} \le \frac{1}{n_k}
\big(\sum^{n_k}_{j=1} \tau_{jk}\big) \sup_{1\le j\le n_k}
\sigma_{jk}=\frac{1}{n_k} \alpha_k\beta_k
\end{equation}
and hence
\begin{equation}
\label{eq-2.13}
n^{-n} \le \prod_{k=1}^m (\alpha_k^{n_k}\beta_k^{n_k}n_k^{-n_k}).
\end{equation}
\eproof

Our next lemma follows from \cite{CP} and Lemma \ref{lemma-2.3}.

\begin{lemma}
\label{lemma-2.4}
There is a universal constant $c>0$ so that if $T_k\in
X_k^*\otimes\ell_2^{n_k}$, $\alpha_k$ is as in Lemma \ref{lemma-2.3} for
$1\le k\le m$ and $T=\oplus^m_{k=1} \alpha_kT_k: X\to \ell_2^n$, then
\begin{equation}
\label{eq-2.14}
c^n \prod^m_{k=1} n_k^{n_k} |T_k(B_{X_k})| \le n^n |T(B_X)|.
\end{equation}
\end{lemma}

\bproof
Let $\{\sigma_{jk}\mid 1\le j\le n_k, 1\le k\le m\}$, $\beta_k$ for
$1\le k\le m$ and $Z$ be defined as in the proof of Lemma
\ref{lemma-2.3} and let $\e_{jk} = \pm 1$ for $1\le j\le n_k$, $1\le
k\le m$. If $\{f_{jk}\mid 1\le j\le n_k, 1\le k\le m\}$ denotes the
canonical basis of $Z^*$ and $B: \ell_2^n\to Z^*$ is the diagonal
operator defined by $(\sigma_{jk})$, then it follows from
\cite[Corollary 1.4(e)]{CP} that there is a universal constant $d>0$ so
that
\begin{eqnarray}
\label{eq-2.15}
d\left(\frac{\prod^m_{k=1}\beta_k^{n_k}}{|B_{Z^*}|}\right)^{\frac{1}{n}}
&=& d\left(\frac{\det(B)}{|B_{Z^*}|}\right)^{\frac{1}{n}} 
\le \hbox{Ave}_\e \big\|\sum^m_{k=1}\sum^{n_k}_{j=1}
\e_{jk}\sigma_{jk} f_{jk}\big\|_{Z^*} \nonumber \\
&=& \big\|\sum^m_{k=1}\sum^{n_k}_{j=1}
\sigma_{kj}f_{jk}\big\|_{Z^*} 
=\big\|\sum^m_{k=1} \beta_ke_k\big\|=1
\end{eqnarray}
and hence
\begin{equation}
\label{eq-2.16a}
d^n\prod^m_{k=1} \beta_k^{n_k} \le |B_{Z^*}|.
\end{equation}
From Lemma \ref{lemma-2.1}, Lemma \ref{lemma-2.3} and (2.14) we
now obtain:
\begin{eqnarray}
\label{eq-2.17a}
|T(B_X)| &=& |B_X||\det(T)|  = |B_X|\prod^m_{k=1} \alpha_k^{n_k} |\det(T_k)|
\nonumber \\
&=& |B_{Z^*}|\Big(\prod^m_{k=1} \alpha_k^{n_k} |\det(T_k)| |B_{X_k}|\Big)
\cdot \Big(\prod^m_{k=1}
|B_{\ell_\infty^{n_k}}|\Big)^{-1}  \nonumber \\
&\ge & 2^{-n} d^n \prod^m_{k=1} n_k^{n_k} |T_k(B_{X_k})| \cdot
\prod^m_{k=1} \alpha_k^{n_k} n_k^{-n_k}\beta_k^{n_k} \nonumber \\
&\ge & 2^{-n} d^nn^{-n} \prod^m_{k=1} n_k^{n_k}|T_k(B_{X_k})|
\end{eqnarray}
which is (\ref{eq-2.14}) with $c=\frac{d}{2}$.
\eproof

\begin{theorem}
\label{thm-2.5}
Let $1\le r,p\le\infty$. There is a universal constant $C$ 
so that if $Y$ is a finite dimensional quotient of a Banach
space $Z$ with $\gl(r,p)$ then
\begin{equation}
\label{eq-2.16}
\max\{\vr(Y,\ell_p),\vr(Y^*,S_{r'})\}\le \vr(Y,\ell_p)\vr(Y^*,S_{r'})\le C\sqrt {r'}\ \gl_{r,p}(Z).
\end{equation}
In particular, if $Y$ is a finite-dimensional quotient of a $r$-convex $(1<r)$
 and $p$-concave Banach lattice $Z$ then
 \[
  \vr(Y,\ell_p)\vr(Y^*,S_{r'})\le C\sqrt {r'}\ K^r(Z)K_p(Z).
\]
\end{theorem}

\bproof
Let $Q: Z\to Y$ be a quotient map, put $n=\dim Y$ and let $S\in
Y^*\otimes\ell_2^n$ be arbitrary. It follows from \cite[Theorem 3.10 i)]
{GJ} that
\begin{eqnarray}
\label{eq-2.17}
n^{\frac12} v_n(S)\vr(Y^*,S_{r'}) &\le& c \sqrt{r'}\ \pi_{r'}(S^*)\le
c\sqrt {r'}\ i_{r'}(Q^*S^*)  \nonumber\\
&\le& c\sqrt {r'}\ \gl_{p',r'}(Z^*)\pi_{p'}(SQ)  \nonumber \\
&\le& c\sqrt {r'}\ \gl_{r,p}(Z)\pi_{p'}(S).
\end{eqnarray}

Since (\ref{eq-2.17}) holds for all $S$ it follows from \cite[Theorem
3.7 (ii)]{GJ} that there is a universal constant $C$ so that
\[
\vr(Y,\ell_p)\vr(Y^*,S_{r'})
\le C\sqrt {r'}\ \gl_{r,p}(Z)
\]
and apply now Theorem \ref{thm-1.2} (vii).
\eproof

\bremark 
We now note that 
\[ \vr(Y,\ell_p)\vr(Y^*,S_{r'})\sim 
\inf\Big(\frac{|V_r|}{|V_p|}\Big)^{\frac{1}{n}}\]
where the infimum is taken over all $n$-dimensional linear quotients $V_r$ 
of $\ell_r$ and $V_p$ of $\ell_p$ so that $V_p\subseteq B_Y\subseteq V_r$. 
Indeed, by definition
\begin{eqnarray*}
\vr(Y,\ell_p) &=& \inf \Big\{ \Big( \frac{|B_Y|}{|V_p|} 
\Big)^{\frac{1}{n}} \mid V_p\subseteq B_Y \Big\} ,\\
\vr(Y^*,S_{r'}) 
&=& 
\inf\Big\{\Big(\frac{|B_{Y^*}|}{|T(B_{S_{r'}})|}\Big)^{\frac{1}{n}} 
\mid S_{r'}\subseteq \ell_{r'},\ T(B_{S_{r'}}) \subseteq B_{Y^*}\Big\}.
\end{eqnarray*}
Now, $[T(B_{S_{r'}})]^o=W(B_{\ell_r})=V_r$ for some linear 
$W:\ell_r\rightarrow\bR^n$. By \cite{S} and \cite{BM} for any 
$n$-dimensional space $X, \ (|B_X||B_{X^*}|)^{\frac{1}{n}}\sim 
\frac{1}{n}$, 
and therefore 
\[ \vr(Y^*,S_{r'})\sim \Big(\frac{|V_r|}{|B_Y|}\Big)^{\frac{1}{n}} \]
where $B_Y\subseteq V_r$, hence the result follows.\eremark

We are now able to prove:

\begin{theorem}
\label{thm-2.6}
Let $E$ be a $m$ dimensional Banach space with 
a normalized $1$-unconditional basis, let $X_k$, $k=1,..,m$, be $n_k$-dimensional
Banach spaces and $X=(\sum_{k=1}^m X_k)_E$.  
Let $1<r\le p < \infty$ or $r=1$ and $p= \infty$. There is a universal 
constant $c_0>0$  and a constant $C(r,p)$ so that for $n=\sum_{k=1}^m n_k$
\begin{eqnarray}
\label{eq-2.18}
\frac{1}{c_0\ p'}\Big(\prod_{k=1}^m \vr(X_k,\ell_p)^{n_k}\Big)^{1/n} \le 
\vr(X,\ell_p) \le  c_{0}C(r,p)\ K^r(E)K_p(E)\big(\prod_{k=1}^m 
\vr(X_k,\ell_p)^{n_k}\big)^{1/n}
\end{eqnarray}
where $C(r,p)= \sqrt{r'}$ for $1<r \le p< \infty$ and $C(1,\infty)=1$
\end{theorem}

\bproof
Let us first prove the left inequality of (\ref{eq-2.18}). By
\cite[Theorem 3.7(ii)]{GJ} in the case $1<p\le\infty$ 
there is a universal constant $A$ 
and operators $T_k\in X^*_k\otimes\ell_2^{n_k}$ with
$\pi_{p'}(T_k)=1$ for $1\le k\le m$ so that
\begin{equation}
\label{eq-2.19}
A^{-1}\vr(X_k,\ell_p)\le n_k|T_k(B_{X_k})|^{1/n_k}\le A\sqrt{p'}\ 
\vr(X_k,\ell_p).
\end{equation}
Let $\alpha_k$, $1\le k\le m$ be chosen as in Lemma \ref{lemma-2.3} and
put $T=\bigoplus_{k=1}^m \alpha_kT_k$. Lemma \ref{lemma-2.2} and Lemma
\ref{lemma-2.3} now give
\begin{equation}
\label{eq-2.20}
\pi_{p'}(T) \le c\sqrt{p'}\ \big\|\sum_{k=1}^m
\alpha_k\pi_{p'}(T_k)e_k^*\big\| = c\sqrt{p'}\ \big\|\sum_{k=1}^m 
\alpha_k e_k^*\big\|=c\sqrt{p'}
\end{equation}
and hence by (2.19), (2.20) and Lemma \ref{lemma-2.4} 
there is a $c'>0$ so that
\begin{eqnarray}
\label{eq-2.21}
c'\big(\prod_{k=1}^m \vr(X_k,\ell_p)^{n_k}\big)^{1/n} & \le & 
c'A\big[\prod_{k=1}^m(n_k^{n_k}|T_k(B_{X_k})|)\big]^{1/n}\nonumber \\
& \le & A n|T(B_X)|^{1/n}\le cA^2p'\ \vr(X,\ell_p)
\end{eqnarray}
which shows the left inequality in (\ref{eq-2.18}).

To prove the right inequality we can by definition find operators
$W_k\in B(\ell_p,X_k)$ for $1\le k\le m$ so that $\|W_k\|=1$ and
\begin{equation}
\label{eq-2.22}
\vr(X_k,\ell_p)=\Big(\frac{|B_{X_k}|}{|W_k(B_{\ell_p})|}\Big)^
{\frac{1}{n_k}}.
\end{equation}
Put for every $1\le k\le m$ \ $Y_k=\ell_p/W_k^{-1} (0)$, let $Q_k$ denote
the quotient map of $\ell_p$ onto $Y_k$, define $V_k\in B(Y_k,X_k)$ so
that $W_k=V_kQ_k$, let $Y=(\sum_{k=1}^m Y_k)_E$ and put 
$V=\oplus^m_{k=1}V_k\in B(Y,X)$.

If $S\in B(\ell_p,Y)$ so that $\|S\|=1$ and
$(\frac{|B_Y|}{|S(B_{\ell_p})|})^{\frac{1}{n}}=\vr(Y,\ell_p)$ then since $\|V\|=1$ we
get using Lemma \ref{lemma-2.1}:
\begin{eqnarray}
\label{eq-2.23}
\vr(X,\ell_p)^n &\le &
\frac{|B_X|}{|VS(B_{\ell_p})|}=(\vr(Y,\ell_p))^n\frac{|B_X|}{|B_Y|}
\frac{|S(B_{\ell_p})|}{|VS(B_{\ell_p})|}\nonumber \\
&= & (\vr(Y,\ell_p))^n|\det V|^{-1} \prod_{k=1}^m
\frac{|B_{X_k}|}{|B_{Y_k}|}\nonumber \\
&=& (\vr(Y,\ell_p))^n|\det V|^{-1} \prod_{k=1}^m
\frac{|B_{X_k}|}{|W_k(B_{\ell_p})|}\cdot \prod_{k=1}^m
\frac{|V_k(B_{Y_k})|}{|B_{Y_k}|}\nonumber \\
&=& (\vr(Y,\ell_p))^n\prod^m_{k=1} \vr(X_k,\ell_p)^{n_k}.
\end{eqnarray}
If $p=\infty$ then $Y_k^*$ is a subspace of an $L_1$-space hence
$\GL(Y_k)=\GL(Y_k^*)=1$. It now follows from Corollary
\ref{cor-1.5} that $Y$ has $\GL$ as well with $\GL(Y)=1$ and
hence by the result of \cite{GMP}
there is a universal constant $C$ so that
\begin{equation}
\label{eq-2.24}
1\le \vr(Y,\ell_\infty)\vr(Y^*,\ell_\infty)\le C\gl(Y)=C.
\end{equation}
If $1<r\le p <\infty$ then it follows from Theorem \ref{thm-1.5} that
$E(\ell_p)$ has $\GL(r,p)$ with $\GL_{r,p}(E(\ell_p))\le K^r(E)K_p(E)$
(this can also easily be obtained from the fact that $E(\ell_p)$ is an
$r$-convex and $p$-concave Banach lattice) and the operator
$Q=\oplus^m_{k=1}Q_k$ is readily seen to be a quotient map of
$E(\ell_p)$ onto $Y$. Hence Theorem \ref{thm-2.5} assures the existence
of a universal constant $C$ so that
\begin{equation}
\label{eq-2.25}
\vr(Y,\ell_p)\le C\sqrt{r'}\ K^r(E)K_p(E).
\end{equation}

Combining (\ref{eq-2.23}) with (\ref{eq-2.24}) and (\ref{eq-2.25}) we
obtain the right inequality of (\ref{eq-2.18}).
\eproof

\vspace{1cm}

\noindent Department of Mathematics,\\ Technion\\ Haifa 32000, Israel,\\
gordon@@tx.technion.ac.il\\

\noindent Department of Mathematics,\\ University of Kiel,\\ Ludewig 
Meyn Strasse 4,\\
D-24098 Kiel, Germany,\\
nms06@@rz.uni-kiel.d400.de\\

\noindent Department of Mathematics
and Computer Science,\\ Odense University,\\
Campusvej 55, DK-5230 Odense M, Denmark \\
and Department of Mathematics,\\ University of Kiel,\\
njn@@imada.ou.dk

\end{document}